\documentclass[10pt,9pt]{article}
\usepackage{amssymb}
\usepackage{makeidx}
\usepackage{amsmath}
\usepackage{mathrsfs}
\usepackage{color}
\usepackage{epsfig}
\newtheorem{example}{Example}[section]
\newtheorem{definition}{Definition}[section]
\newtheorem{theorem}{Theorem}[section]

\newtheorem{lemma}{Lemma}[section]

\begin{document}
\title{Compactly Supported Cohomology Groups of Smooth Toric Surfaces}

\author{Ma\l gorzata Aneta Marciniak}
\maketitle

{\it Abstract:} This article uses homological methods for evaluating compactly supported cohomology groups of noncompact smooth toric surfaces.

\bigskip


{\it Kyewords:} Compactly supported cohomology groups, toric varieties

\bigskip

\section{Introduction}

Compactly supported cohomology groups play a key role in complex analysis. Vanishing properties are related to solvability of $\overline{\partial}$ problem as well as analytic continuation problems (\cite{dw-1} and \cite{dw-2}). This article uses the additive property for finding the compactly supported cohomology groups of smooth toric surfaces.
The sheaf of germs of holomorphic functions on $X$ is denoted as $\mathscr{O}_X$ or simply $\mathscr{O}$, if there is no confusion which space is considered.

Compactifications and partial compactifications of toric surfaces play a key role in the procedure of computing compactly supported cohomology groups. Let $X$ be a noncompact toric surface associated with the fan $\Sigma$. Let us denote a toric compactification of $X$ as $\widetilde X $ with the fan $\widetilde\Sigma$. The approach to the problem varies depending on convexity properties of the connected components of the (open) set $supp(\widetilde{\Sigma})\setminus supp(\Sigma)$. The results are described in three theorems:
\begin{enumerate}
\item{} Exactly one component of $supp(\widetilde{\Sigma})\setminus supp(\Sigma)$ is concave (Theorem \ref{thm:strictly_convex}). This condition is equivalent to the fan $\Sigma$ being a subfan of a strictly convex fan.
\item{} Exactly one component of $supp(\widetilde{\Sigma})\setminus supp(\Sigma)$ is a half-plane (Theorem \ref{thm:half_plane}). This condition is equivalent to $\Sigma$ being a subfan of a fan that covers a half plane.
\item{} Exactly two components of $supp(\widetilde{\Sigma})\setminus supp(\Sigma)$ are half-planes (Example \ref{ex:line}). This condition implies that the toric surface is $\mathbb{P}^1\times \mathbb{C}^1$.
\item{} All components of $supp(\widetilde{\Sigma})\setminus supp(\Sigma)$ are strictly convex (Theorem \ref{thm:entire plane}). It can be seen as the fan $\Sigma$ spans the entire plane.
\end{enumerate}
In this terminology a component is strictly convex if it does not contain a line.

An introduction to the theory of toric varieties can be found in \cite{bonavero} or in \cite{dw}. Toric surfaces are presented shortly in \cite{marciniak}.

 %Similar results can be formulated for other function sheafs.

\section{Compactly Supported Cohomology Groups}
This section contains a short overview of definitions and main properties of the compactly supported cohomology groups.

\subsection{Definition and Examples}

\begin{definition}{\bf (Compactly supported Dolbeault cohomology groups)}
Compactly supported Dolbeaut cohomology groups of the domain $D$ are the complex vector spaces:
\[
\mathfrak{H}_{c}^{p,q}(D)=\frac{\{\overline{\partial}\mbox{-closed forms with compact support of bidegree $(p,q)$ in $D$}\}}{\{\overline{\partial}\mbox{-exact forms with compact support of bidegree $(p,q)$ in $D$}\}}.
\]
\end{definition}
\bigskip
The following theorem shows the relationship between compactly supported cohomology groups and compactly supported Dolbeault cohomology groups.

\begin{theorem}{\bf (Dolbeault's Theorem, \cite{do})}
If $D$ is an open domain in the space of $n$ complex variables, $\mathscr{O}$ is the sheaf of germs of holomorphic functions on $D$, and $\mathfrak{H}_{c}^{p,q}(D)$ is compactly supported Dolbeault cohomology group of bidegree $(p,q$) for $D$. Then $H_{c}^q(D,\mathscr{O})=\mathfrak{H}_{c}^{0,q}(D)$.\,\rule[-.2ex]{.6ex}{1.8ex}
\end{theorem}

An alternate definition can be found in \cite{bredon}.
In particular, if $X$ is a compact manifold, then
\[
H^i_c(X,\mathscr{O})=H^i(X,\mathscr{O})=
\begin{cases}
0          &\text{if $i\neq 0$}\\
\mathbb{C} &\text{$i=0$,}
\end{cases}
\]
i.e., compactly supported groups are equal to the usual cohomology groups for compact spaces.
It is a well known fact that the groups $H_c^i(\mathbb{C}^n,\mathscr{O})$ are trivial for all $n\geq 2$ and $i\geq 0$ (\cite{krantz}).

\subsection{Properties}

Compactly supported cohomology groups have all properties of a cohomology theory.
The ``additive" property, broadly used in the further part of the research, requires the notion of the inverse image of a sheaf.

\begin{definition}\textbf{(Inverse image)}
Let $f:A\rightarrow B$ be a map and let $\mathscr{G}$ be a sheaf on $B$ with canonical projection $\pi:\mathscr{G}\rightarrow B$. The inverse image sheaf $f^{*}\mathscr{G}$ is defined as
\[
f^{*}\mathscr{G}=\{(a,g)\in A\times \mathscr{G}:f(a)=\pi(g)\}.
\]
\end{definition}
In particular, if $f$ is a closed embedding, the following theorem holds.

\begin{theorem} \textbf{(\cite{iversen} III.7.6)}
Let $i:Y\rightarrow X$ be a closed embedding, then the following sequence
\[
\ldots\rightarrow H_c^q(X\setminus Y, \mathscr{F}) \rightarrow
H_c^q(X,\mathscr{F}) \rightarrow
H_c^q(Y, i^{*} \mathscr{F})\rightarrow
H_c^{q+1}(X\setminus Y, \mathscr{F})\rightarrow \ldots,
\]
is exact. \,\rule[-.2ex]{.6ex}{1.8ex}
\end{theorem}

The following exact sequences are obtained for each $n$ separately.
\begin{theorem}\textbf{(The K\"unneth Formula, \cite{bredon} II. Theorem 15.2)}
If $X$ and $Y$ are locally compact Hausdorff spaces, with the sheaves $\mathscr{F}$ and $\mathscr{G}$ respectively and $\mathscr{F}*\mathscr{G}=0$, then the sequence
\[
0\rightarrow \bigoplus_{p+q=n} H^p_c(X,\mathscr{F})\otimes H^q_c(Y,\mathscr{G})\rightarrow H^n_c(X\times Y,\mathscr{F}\otimes \mathscr{G})\rightarrow\bigoplus_{p+q=n+1} H^p_c(X,\mathscr{F})* H^q_c(Y,\mathscr{G})\rightarrow 0
\]
is exact.
\end{theorem}

\section{Compactly Supported Cohomology Groups of Toric Surfaces}
Main results are written in the form of three theorems, depending on if the support of the fan of the toric surface spans less than a half plane, a half-plane, or the entire plane. The case when the fan spans a line requires entirely different approach and is described in Example \ref{ex:line}. All previous cases use a compactification or a partial compatification of the toric surface. If the fan $\Sigma$ is a subfan of a fan that is strictly convex or covers a half plane then a partial compactification will be sufficient. For all cases the notation is the same.

Let $X$ be a smooth noncompact toric variety associated with the fan $\Sigma$ and let its toric compatification $\widetilde{X} $ be associated with the fan $\widetilde{\Sigma} $. Clearly, $\Sigma$ is a subfan of $\widetilde{\Sigma}$ and we can consider components of $supp(\overline{\Sigma})\setminus supp(\Sigma)$. Let $C_1,\ldots,C_n$ be strictly convex components of $supp(\overline{\Sigma})\setminus supp(\Sigma)$ and let $C_0$ be the concave component or an (open) half plane (the component $C_0$ can be omitted since it does not  appear in the cohomology formula).
Then $\widetilde{X}\setminus X$ consists of connected components $Y_1,\ldots, Y_n$ that are defined by the components $C_1,\ldots,C_n$. The embeddings of $Y_j$ into $\widetilde{X} $ will be denoted as $i_j:Y_j\rightarrow \widetilde{X}$ and the (strictly) convex connected components of $\text{supp}(\widetilde{\Sigma})\setminus \text{supp}(\Sigma)$ are spanned by the pairs of vectors $v_j$ and $w_j$ (with the positive orientation of $\mathbb{R}^2$) and have the (singularity) type $(p_j,q_j)$ as described for toric surfaces in \cite{ful} section 2.6.

\subsection{The Fan Spans Less than a Half-Plane}

The following theorem describes compactly supported cohomology groups of a smooth toric surface $X$ with a fan $\Sigma$ that spans less than a half plane. In other words, $\Sigma$ can be seen as a subfan of a strictly convex fan.

\begin{theorem}\label{thm:strictly_convex}
Let $X$ be a smooth toric surface which fan $\Sigma$ is a subfan of a strictly convex fan. Then $H^0_c(X,\mathscr{O})=H^2_c(X,\mathscr{O})=0$ and
\[
H^1_c(X,\mathscr{O})=\bigoplus_{j=1}^{n}H^0_c(Y_j,i_j^{*}\mathscr{O})=\bigoplus_{j=1}^{n}\{\sum_{(s,t)>(0,0)} a_{st}z_j^sw_j^t: {p_j}t\geq {q_j} s\},
\]
where $Y_j$ are varieties defined by (strictly) convex connected components of $\mathbb{R}^2\setminus \text{supp}(\Sigma)$ and $i_j$ are (closed) embeddings $i_j:Y_j\rightarrow X$ and the series converge.
\end{theorem}

Proof. Since $\Sigma$ is a subfan of a strictly convex fan, $X$ can be treated as an open submanifold of a smooth toric surface $\widetilde{X}$ with a strictly convex fan $\widetilde{\Sigma}$. Then $X=\widetilde{X}\setminus Y$, where $Y$ is a subvariety of $\widetilde{X}$, and $Y=Y_1\cup \ldots \cup Y_n$ is its decomposition into disjoint, compact and connected subvarieties. Note that the additive property provides the following exact sequence.
\[
0\rightarrow H^0_c(X,\mathscr{O})\rightarrow H^0_c(\widetilde{X},\mathscr{O})\rightarrow H^0_c(Y,i^{*}\mathscr{O})\rightarrow H^1_c(X,\mathscr{O})\rightarrow
H^1_c(\widetilde{X}, \mathscr{O})\rightarrow H^1_c(Y,i^{*}\mathscr{O})\rightarrow  H^2_c(X,\mathscr{O})\rightarrow H^2_c(\widetilde{X},\mathscr{O})\rightarrow 0.
\]
Then $H^0_c(X,\mathscr{O})= H^0_c(\widetilde{X},\mathscr{O})=0$ because $X$ and $\widetilde{X}$ are noncompact; $H^0_c(Y,i^{*}\mathscr{O})=H^0_c(Y,\mathscr{O}_{Y})=\bigoplus_{j=1}^{n}H^0_c(Y_j,\mathscr{O}_{Y_j})$ because $Y_j$ are compact and disjoint; $H^1_c(\widetilde{X}, \mathscr{O})=0$ from Theorem 4.2 in \cite{marciniak}; and $H^2_c(\widetilde{X},\mathscr{O})=0$ because of dimensional reasons. Thus the sequence simplifies to the following two sequences:
\[
0\rightarrow H^0_c(Y,i^{*}\mathscr{O})\rightarrow H^1_c(X,\mathscr{O})\rightarrow 0
\]
and
\[
0\rightarrow H^1_c(Y,i^{*}\mathscr{O})\rightarrow  H^2_c(X,\mathscr{O})\rightarrow H^2_c(\widetilde{X},\mathscr{O})\rightarrow 0.
\]
Then
\[
 H^1_c(X,\mathscr{O})= H^0_c(Y,i^{*}\mathscr{O})=\bigoplus_{j=1}^{n} H^0_c(Y_j,i^{*}\mathscr{O}),
\]
since $Y_j$ are disjoint and compact components of $Y$ and
\[
H^2_c(X,\mathscr{O})= H^1_c(Y,i^{*}\mathscr{O})=\bigoplus_{j=1}^{n}H^1_c(Y_j,i^{*}\mathscr{O})=0,
\]
since each $Y_j$ is compact (as the sum of projective curves) in $\widetilde{X}$. \,\rule[-.2ex]{.6ex}{1.8ex}

\bigskip

Note that the groups $ H^0_c(Y_j,i^{*}\mathscr{O})$ can be found explicitly in terms of $p_j$ and $q_j$. If $(z_j,w_j)$ are local coordinates then
\[
H^0_c(Y_j,i^{*}\mathscr{O})=\{\sum_{(s,t)> (0,0)} a_{st}z_j^sw_j^t: {p_j}t\geq {q_j} s\},
\]
where the series $\displaystyle{\sum_{(s,t)>(0,0)} a_{st}z_j^sw_j^t}$ converges in the coordinates $(z_j,w_j)$ in some neighborhood of $Y_j$ in $\widetilde{X}$. The details can be found in Section 2.2 of \cite{marciniak2}.

\subsection{The Fan Spans a Half-Plane}

The following lemma evaluates the compactly supported cohomology groups for those toric surfaces which fans have the supports that covers a half-plane. We will need this results for the next theorem.

\begin{lemma}\label{lemma:half-plane}
Let $X$ be a smooth toric surface which fan $\Sigma$ such that $\text{supp}(\Sigma)$ is a half-plane. Then $H^0_c(X,\mathscr{O})=H^2_c(X,\mathscr{O})=0$ and $H^1_c(X,\mathscr{O})=\{\sum_{s> 0} a_{s}z^s\}$, where the series converges in a neighborhood of $0$.
\end{lemma}
Proof. Note that that $X$ can be represented as $X=\widetilde{X}\setminus \mathbb{P}^1$, where $\widetilde{X}$ is a smooth compact toric surface. Then $X$ admits the following exact sequence:
\[
0\rightarrow H^0_c(X,\mathscr{O})\rightarrow H^0_c(\widetilde{X},\mathscr{O})\rightarrow H^0_c(\mathbb{P}^1,i^{*}\mathscr{O})\rightarrow H^1_c(X,\mathscr{O})\rightarrow
H^1_c(\widetilde{X}, \mathscr{O})\rightarrow H^1_c(\mathbb{P}^1,i^{*}\mathscr{O})\rightarrow  H^2_c(X,\mathscr{O})\rightarrow H^2_c(\widetilde{X},\mathscr{O})\rightarrow 0.
\]
Note that $H^0_c(X,\mathscr{O})=0$ since $X$ is noncompact and $H^0_c(\widetilde{X},\mathscr{O})=\mathbb{C}$ since $\widetilde X$ is compact. Moreover $H^1_c(\widetilde{X}, \mathscr{O})= H^1_c(\mathbb{P}^1,i^{*}\mathscr{O})=  H^2_c(X,\mathscr{O})=H^2_c(\widetilde{X},\mathscr{O})= 0$ and the sequence simplifies to:
\[
0\rightarrow \mathbb{C}\rightarrow H^0_c(\mathbb{P}^1,i^{*}\mathscr{O})\rightarrow H^1_c(X,\mathscr{O})\rightarrow 0.
\]
Then the group $H^1_c(X,\mathscr{O})$ is a quotient of $H^0_c(\mathbb{P}^1,i^{*}\mathscr{O})$ and $\mathbb{C}$. Since the embedding of $\mathbb{P}^1$ into $\widetilde X$ is flat we obtain that $H^0_c(\mathbb{P}^1,i^{*}\mathscr{O})=\{\sum_{s\geq 0} a_{s}z^s\}$ and $H^1_c(X,\mathscr{O})=\{\sum_{s> 0} a_{s}z^s\}$, where the series converge in a neighborhood of $0$.
\,\rule[-.2ex]{.6ex}{1.8ex}

Note that $H^1_c(X,\mathscr{O})=H^1_c(\mathbb{C}^1,\mathscr{O})$.
\bigskip

Now the following theorem can be formulated.

\begin{theorem}\label{thm:half_plane}
Let $X$ be a smooth toric surface which fan $\Sigma$ is a subfan of a fan $\widetilde{\Sigma}$ which support is a half-plane. Then $H^0_c(X,\mathscr{O})=H^2_c(X,\mathscr{O})=0$ and
\[
H^1_c(X,\mathscr{O})=H^1_c(\mathbb{C}^1, \mathscr{O})\oplus \bigoplus_{j=1}^{n}H^0_c(Y_j,i^{*}_j\mathscr O)=\{\displaystyle{\sum_{s>0} a_{s}z_0^s }\}\oplus \bigoplus_{j=1}^{n}\{\sum_{(s,t)> (0,0)} a_{st}z^sw^t: {p_j}t\geq {q_j} s\},
\]
where all series converge and the (strictly) convex connected components of $\text{supp}\widetilde{\Sigma}\setminus\text{supp}(\Sigma)$ define the varieties $Y_i$.
\end{theorem}

Proof. Let $\widetilde{X}\setminus X= Y$, where $\widetilde{X}$ is a smooth compact toric variety and $Y=Y_1\cup\ldots\cup Y_n$ is the decomposition of $Y$ into compact and connected components. Then
\[
0\rightarrow H^0_c(X,\mathscr{O})\rightarrow H^0_c(\widetilde{X},\mathscr{O})\rightarrow H^0_c(Y,i^{*}\mathscr{O})\rightarrow H^1_c(X,\mathscr{O})\rightarrow
H^1_c(\widetilde{X}, \mathscr{O})\rightarrow H^1_c(Y,i^{*}\mathscr{O})\rightarrow  H^2_c(X,\mathscr{O})\rightarrow H^2_c(\widetilde{X},\mathscr{O})\rightarrow 0.
\]
Then: $H^0_c(X,\mathscr{O})= H^0_c(\widetilde{X},\mathscr{O})=0$ because $X$ and $\widetilde{X}$ are noncompact; $H^0_c(Y,i^{*}\mathscr{O})=\bigoplus_{j=1}^{n}H^0_c(Y_j,i^{*}_j\mathscr O)$ because $Y_j$ are compact and disjoint; $H^1_c(\widetilde{X}, \mathscr{O})=\{\sum_{s> 0} a_{s}z^s\}$ from Lemma \ref{lemma:half-plane}; $H^2_c(\widetilde{X},\mathscr{O})=0$ from Lemma \ref{lemma:half-plane}; and $H^1_c(Y,i^{*}\mathscr{O})=0$ as in the proof of the previous theorem. Thus the sequence simplifies to:
\[
0\rightarrow H^0_c(Y,i^{*}\mathscr{O})\rightarrow H^1_c(X,\mathscr{O})\rightarrow
H^1_c(\widetilde{X}, \mathscr{O})\rightarrow 0\rightarrow  H^2_c(X,\mathscr{O})\rightarrow 0.
\]
Then $H^2_c(X,\mathscr{O})=0$ and
\[
0\rightarrow \bigoplus_{j=1}^{n}H^0_c(Y_j,i^{*}_j\mathscr O)\rightarrow H^1_c(X,\mathscr{O})\rightarrow
H^1_c(\widetilde{X}, \mathscr{O})\rightarrow 0.
\]
The sequence does not split but we still obtain the following result
\[
H^1_c(X,\mathscr{O})=H^1_c(\mathbb{C}^1, \mathscr{O})\oplus \bigoplus_{j=1}^{n}H^0_c(Y_j,i^{*}_j\mathscr O)=\{\sum_{s>0} a_{s}z_0^s\}\oplus \bigoplus_{j=1}^{n}\{\sum_{(s,t)> (0,0)} a_{st}z_j^sw_j^t: {p_j}t\geq {q_j} s\}.
\]
\,\rule[-.2ex]{.6ex}{1.8ex}

The details of the last step are presented is Section 2.2 of \cite{marciniak2}.

\subsection{The Fan Spans a Line}

If the support of the fan $\Sigma$ is a line, then the toric surface is simply $\mathbb{P}^{1}\times \mathbb{C}^{*}$ and its cohomology groups can be computed from the K\"unneth formula.

\begin{example}\label{ex:line}
Let us find $H^i_c(\mathbb{P}^1\times \mathbb{C}^*,\mathscr{O})$ for $i=0,1,2$ using the K\"unneth formula for products. Recall that $H^0_c(\mathbb{P}^1,\mathscr{O})=\mathbb{C}$, $H^1_c(\mathbb{P}^1,\mathscr{O})=0$, $H^0_c(\mathbb{C}^*,\mathscr{O})=0$ and $H^1_c(\mathbb{C}^*,\mathscr{O})=\{\sum_{s> 0} a_{s}z^s, a_s\in \mathbb{C}\}\oplus\{\sum_{s> 0} a_{s}\frac{1}{z^s}, a_s\in \mathbb{C}\}$. Then the K\"unneth Formula for $\mathbb{P}^1\times \mathbb{C}^*$ gives the following for the first cohomology group of $\mathbb{P}^1\times \mathbb{C}^*$:
\[
0\rightarrow \bigoplus_{p+q=1} H^p_c(\mathbb{C}^1,\mathscr{O})\otimes H^q_c(\mathbb{P}^1,\mathscr{O})\rightarrow H^1_c(\mathbb{C}^*\times \mathbb{P}^1,\mathscr{O})\rightarrow\bigoplus_{p+q=2} H^p_c(\mathbb{C}^*,\mathscr{O})* H^q_c(\mathbb{P}^1,\mathscr{O})\rightarrow 0,
\]
which converts to:
\[
0\rightarrow H^1_c(\mathbb{C}^*,\mathscr{O})\bigoplus H^0_c(\mathbb{P}^1, \mathscr{O})\rightarrow H^1_c(\mathbb{C}^*\times \mathbb{P}^1,\mathscr{O})\rightarrow H^1_c(\mathbb{C}^*,\mathscr{O})* H^1_c(\mathbb{P}^1,\mathscr{O})\rightarrow 0,
\]
and proves that $H^1_c(\mathbb{C}^*\times \mathbb{P}^1,\mathscr{O})=H^1_c(\mathbb{C}^*,\mathscr{O})=\{\sum_{s> 0} a_{s}z^s \}\oplus\{\sum_{s> 0} a_{s}\frac{1}{z^s}, a_s\in \mathbb{C}\}$. Similarly for $H^2_c(\mathbb{C}^*\times \mathbb{P}^1,\mathscr{O})$:
\[
0\rightarrow H^1_c(\mathbb{C}^*,\mathscr{O})\otimes H^1_c(\mathbb{P}^1,\mathscr{O}) \rightarrow H^2_c(\mathbb{C}^*\times \mathbb{P}^1,\mathscr{O})\rightarrow 0,
\]
which implies $H^2_c(\mathbb{C}^*\times \mathbb{P}^1,\mathscr{O})= 0$.
\end{example}

\subsection{The Fan Spans the Entire Plane}

If the support of the fan $\Sigma$ is not a subset of a half-plane then it spans the entire plane. In this case all components of $\text{supp}\widetilde{\Sigma}\setminus \text{supp}\Sigma$ are strictly convex.

\begin{theorem}\label{thm:entire plane}
Let $X$ be a noncompact toric surface which fan $\Sigma\subset \mathbb{R}^n$ has the support that is not a subset of a half-plane. Then $H^0_c(X,\mathscr{O})=H^2_c(X,\mathscr{O})=0$ and
\[
H^1_c(X,\mathscr{O})=\bigoplus_{j=1}^{n}H^0_c(Y_j,i^{*}_j\mathscr O)\bigoplus_{j=1}^{n}H^0_c(Y_j,i_j^{*}\mathscr{O})=\bigoplus_{j=1}^{n}\{\sum_{(s,t)>(0,0)} a_{st}z_j^sw_j^t: {p_j}t\geq {q_j} s\},
\]
where strictly convex connected components of $\text{supp}\widetilde{\Sigma}\setminus\text{supp}(\Sigma)$ define the varieties $Y_j$.
\end{theorem}

Proof. Let $\widetilde{X}\setminus X =Y$, where $\widetilde{X}$ is a smooth compact toric variety and $Y=Y_1\cup\ldots\cup Y_n$ is the decomposition of $Y$ into compact and connected components. Then the additive property has the following form:
\[
0\rightarrow H^0_c(X,\mathscr{O})\rightarrow H^0_c(\widetilde{X},\mathscr{O})\rightarrow H^0_c(Y,i^{*}\mathscr{O})\rightarrow H^1_c(X,\mathscr{O})\rightarrow
H^1_c(\widetilde{X}, \mathscr{O})\rightarrow H^1_c(Y,i^{*}\mathscr{O})\rightarrow  H^2_c(X,\mathscr{O})\rightarrow H^2_c(\widetilde{X},\mathscr{O})\rightarrow 0.
\]
Recall that $H^0_c(\widetilde{X},\mathscr{O})=\mathbb{C}$ and $H^1_c(\widetilde{X},\mathscr{O})=H^2_c(\widetilde{X},\mathscr{O})=0$ since $\widetilde{X}$ is compact. Moreover, since $X$ is noncompact, $H^0_c(X,\mathscr{O})=0$. Thus the sequence can be written as two sequences:
\[
0\rightarrow \mathbb{C}\rightarrow H^0_c(Y,i^{*}\mathscr{O})\rightarrow H^1_c(X,\mathscr{O})\rightarrow 0
\]
and
\[
0\rightarrow H^1_c(Y,i^{*}\mathscr{O})\rightarrow  H^2_c(X,\mathscr{O})\rightarrow 0.
\]
Then $H^2_c(X,\mathscr{O})=H^1_c(Y,i^{*}\mathscr{O})$ and
\[
H^1_c(X,\mathscr{O})=H^0_c(Y,i^{*}\mathscr{O})/\mathbb{C}=\bigoplus_{j=1}^{n}H^0_c(Y_j,i_j^{*}\mathscr{O})/\mathbb{C}=\bigoplus_{j=1}^{n}\{\sum_{(s,t)>(0,0)} a_{st}z_j^sw_j^t: {p_j}t\geq {q_j} s\}
\]
\,\rule[-.2ex]{.6ex}{1.8ex}

\section{Further Research}

If $\Sigma\subset \mathbb{R}^m$ is a fan associated with a noncompact toric variety $X$ then the compactly supported cohomology groups as well as the obstacles for solving extension problems lie in the connected components of $\mathbb{R}^m\setminus \text{supp}(\Sigma)$. In terms of the geometry (or topology) of $X$ those components describe the ends of $X$. The notion of an end was formally introduced by Freudenthal in \cite{fr} and is commonly used by analysts, for example in \cite{ghu} or \cite{marciniak}. Asking which properties of ends in higher dimensions have impact on the cohomology groups yields to interesting questions that involve compactifications of toric varieties and tropical geometry.


\begin{thebibliography}{BBBC}

\bibitem{bonavero} L.~Bonavero and M.~Brion (Eds.),
{\it Geometry of toric varieties. Lectures from the Summer
School held in Grenoble},
June 19--July 7, 2000. S\'eminaires et Congres, 6.
Soci\'et\'e Math\'ematique de France, Paris (2002).

\bibitem{bredon} G. E.~Bredon, {\it Sheaf Theory}, 2 ed., Springer, (1997).

\bibitem{do} P.~Dolbeault, {\it Sur la cohomologie des vari\' et\' es analytiques complexes}
C. R. Acad. Sci. Paris, 236 (1953) 175 - 177.

\bibitem{dw} R.~Dwilewicz,
{\it An analytic point of view at toric varieties},
Serdica Math. J. {\bf 33} (2007), 163 - 240.

\bibitem{dw-1} R.~Dwilewicz,
{\it Additive Riemann-Hilbert problem
in line bundles over $\mathbb{CP}^1$},
Canadian Math. Bull. {\bf 49} (2006), 72 - 81.


\bibitem{dw-2} R.~Dwilewicz,
{\it Holomorphic extensions in complex fiber bundles},
J.~Math.~Analysis and Appl. {\bf 322} (2006),

\bibitem{ghu} B.~Gilligan, A.T.~Huckleberry, {\it Complex homogeneous manifolds with two ends},
Michigan Math. J. (1981).

\bibitem{fr} H.~Freudenthal, {\it \"Uber die Enden topologischer
R\"aume und Gruppen}, Math.~Z. {\bf 33} (1931), 692-713.

\bibitem{ful} W.~Fulton, {\it Introduction to Toric Varieties}, Princeton Univ.~Press, NJ, (1993).

\bibitem{iversen} B.~Iversen, {\it Cohomology of Sheaves},
 Springer-Verlag, (1986).

\bibitem{krantz} S.G.~Krantz, {\it Function
Theory of Several Complex Variables}, AMS Chelsea
Publishing, (2001).

\bibitem{marciniak} M.~A.~Marciniak, {\it Holomorphic Extensions in Toric Surfaces}, Journal of Geometric Analysis 2011, 1-23.


\bibitem{marciniak2} M.~A.~Marciniak, {\it Using the Additive Property of Compactly Supported Cohomology Groups}, 
http://arxiv.org/abs/1209.1035


\end{thebibliography}
\end{document}